\documentclass[11pt]{article}
\usepackage{graphicx}
\usepackage{amsmath}
\usepackage{amsfonts}
\usepackage{epsfig}
\usepackage{setspace}

\newcommand{\be}{\begin{equation}}
\newcommand{\ee}{\end{equation}}
\newcommand{\beqn}{\begin{eqnarray}}
\newcommand{\eeqn}{\end{eqnarray}}
\newcommand{\beqns}{\begin{eqnarray*}}
\newcommand{\eeqns}{\end{eqnarray*}}

\newcommand{\card}{\mbox{card}}

\newcommand{\Cov}{\mbox{Cov}\ }

\newcommand{\EE}{\ensuremath{{\mathbb E}}}
\newcommand{\II}{\ensuremath{{\mathbb I}}}

\newcommand{\fr}[1]{(\ref{#1})}

\newcommand{\om}{\omega}

\newcommand{\Om}{\Omega}
\newcommand{\Te}{\Theta}

\newcommand{\bom}{\mbox{\mathversion{bold}$\om$}}

\newtheorem{lemma}{Lemma}
\newtheorem{theorem}{Theorem}

\newtheorem{remark}{Remark}

\begin{document}

\title{\Large{\bf Anisotropic functional deconvolution with long-memory noise: the case of a multi-parameter fractional Wiener sheet }}

\author{
\large{ Rida Benhaddou}  \footnote{E-mail address: Benhaddo@ohio.edu} \ and\   Qing Liu
  \\ \\
Department of Mathematics, Ohio University, Athens, OH 45701} 
\date{}

\doublespacing
\maketitle
\begin{abstract}
We look into the minimax results for the anisotropic two-dimensional functional deconvolution model with the two-parameter fractional Gaussian noise. We derive the lower 
bounds for the $L^p$-risk, $1 \leq p < \infty$, and taking advantage of the Riesz poly-potential, we apply a wavelet-vaguelette expansion to de-correlate the anisotropic fractional Gaussian noise. We construct an adaptive wavelet hard-thresholding estimator that attains asymptotically quasi-optimal convergence rates in a wide range of Besov balls. Such convergence rates depend on a delicate balance between the parameters of the Besov balls, the degree of ill-posedness of the convolution operator and the parameters of the fractional Gaussian noise. A limited simulations study confirms theoretical claims of the paper. The proposed approach is extended to the general $r$-dimensional case, with $r> 2$, and the corresponding convergence rates do not suffer from the curse of dimensionality. \\


{\bf AMS (2000) Subject Classification: 62G05, 62G20, 62G08 }
 \end{abstract} 

\section{Introduction.}
Consider the problem of estimating a periodic two-dimensional function, $f(\cdot, \cdot)\in L^2(U)$, based on noisy convolutions that are described by the equations
\be \label{conveq}
dY(t, x)=q(t, x)dtdx + \varepsilon^{\overline{\alpha}}dB^{\bf{\alpha}}(t, x), \ q(t, x)=\int^1_0f(s, x)g(t-s, x)ds.
\ee
Here, $U=[0, 1]^2$, $g(\cdot , \cdot)$ is the convolution kernel and it is assumed to be known, $B^{{\alpha}}(t, x)$ is an anisotropic two-dimensional fractional Brownian sheet (fBs), $\overline{\alpha}={(\alpha_1+ \alpha_2)}/{2}$, and $\alpha_i=2-2H_i\in(0, 1]$, $i=1,2$, are the parameters of the long-memory in the direction of $t$ and $x$, respectively. A two-dimensional fractional Brownian sheet is defined by the formula
\be \label{fb-sb}
B^{{\alpha}}(t, x)=\frac{1}{C_H}\int^{t}_0\int^{x}_0Q(s, y)dW(s, y),
\ee
where $W=\left\{ W(t, x), (t, x) \in U \right\}$ is a two-dimensional standard Brownian sheet, $C_H$ is some explicit constant and
\be \label{fb-Q}
Q(s, y)=\left[(t-s)_{+}\right]^{H_1-1/2}\left[(x-y)_{+}\right]^{H_2-1/2},
\ee
with $Y_{+}=\max(0, Y)$. In addition, the two-parameter fBs is characterized by a covariance function of the form 
\be \label{cov}
\Cov(B^{{\alpha}}(t_1, s_1), B^{{\alpha}}(t_2, s_2)) =c\left[|t_1|^{2H_1}+|t_2|^{2H_1}-|t_1-t_2|^{2H_1}\right] \left[|s_1|^{2H_2}+|s_2|^{2H_2}-|s_1-s_2|^{2H_2}\right],
\ee 
for some $c\in R$, $(t_i, s_i)\in U$, $i=1, 2$, where $H_i=1-\frac{\alpha_i}{2}$ are the Hurst parameters in the direction of $t$ and $s$, respectively (see Kamont~(1996)). \\
The discrete, and the more realistic, version of model \fr{conveq} is given by
\be \label{disconveq}
Y(t_i, x_l)=q(t_i, x_l) + \sigma \xi_{il},\ (i, l)=1, 2, \cdots, N,
\ee
where $t_i=\frac{i}{N}$, $x_l=\frac{l}{N}$, $\sigma$ is a positive variance constant, and  $\xi_{il}$ are zero-mean second-order stationary Gaussian random variables satisfying the auto-covariance structure
\be \label{auto-cov}
\gamma(h_1, 0)\asymp h_1^{-\alpha_1},\ \  \gamma(0, h_2)\asymp h_2^{-\alpha_2} \ and\ \ \gamma(h_1, h_2)\asymp h_1^{-\alpha_1}h_2^{-\alpha_2}.
\ee
\noindent
{\bf Assumption A.1.  }\label{A1} The error structure $\{\xi_{il}: i, l \in Z\}$ in model \fr{disconveq} is a zero-mean second-order stationary process satisfying
\be
X_N\xrightarrow{D} B^{{\alpha}},
\ee 
where
\be
X_N(t, x)=\frac{1}{N^{2-\alpha_1/2 -\alpha_2/2}}\sum^{\lceil Nt \rceil}_{i=1}\sum^{\lceil Nx \rceil}_{l=1}\xi_{il}, \ (t, x)\in U^2,
\ee
and $B^{{\alpha}}$ is a two-parameter fractional Brownian sheet.\\
Assumption {\bf A.1} is valid under auto-covariance structure \fr{auto-cov} and it appeared in Adu and Richardson~(2018). It guarantees that, with the calibration $\varepsilon^{{\alpha_1 + \alpha_2}}\asymp n^{\alpha_1/2 + \alpha_2/2}$, models \fr{conveq} and \fr{disconveq} are asymptotically equivalent. Therefore, for sufficiently large sample size $n=N^2$, model \fr{conveq} can be used to approximate model \fr{disconveq}.

Deconvolution model has been the subject of a great deal of papers since late 1980s, but the most significant contribution was that of Donoho~(1995) who was the first to devise a wavelet solution to the problem. Other attempts include,  Abramovich and  Silverman (1998), Walter  and  Shen (1999), Johnstone et al.~(2004), Donoho  and  Raimondo (2004), among others.  In the case of functional deconvolution model with $f(t, x)\equiv f(t)$, Pensky and sapatinas~(2009, 2010, 2011) pioneered into the formulation and further development of the problem. 
 
Functional deconvolution problem of type \fr{conveq} with $\alpha_1= \alpha_2=1$ corresponds to the white noise case and it was investigated under the $L^2$-risk in Benhaddou et al.~(2013), and under $L^p$-risk, $1\leq p < \infty$, in Benhaddou~(2017), where they constructed an adaptive hard-thresholding wavelet estimator, and showed that it is asymptotically quasi-optimal within a logarithmic factor of $\varepsilon$ over a wide range of Besov balls of mixed regularity. This model is motivated by experiments 
in which one needs to recover a two-dimensional function using  observations  of its convolutions 
along profiles $x=x_l$. This situation occurs, for example, 
in seismic inversions (see Robinson~(1999)).  In these articles, it is assumed that the error terms are white noise processes or i.i.d noise. However, empirical evidence has shown that, even at large lags, the correlation structure in the errors takes the power-like form \fr{auto-cov}. This phenomenon is referred to as long-memory (LM) or long-range dependence (LRD). The presence of long memory in oceanic seismic data was pointed out in Wood et al.~(2014) for instance, and it may be due for example to sea floor temperature variations.

Long-memory has been investigated quite considerably  in the standard deconvolution model. One can list; Wang~(1996, 1997), Wishart~(2013), Benhaddou et al.~(2014) and Kulik et al.~(2015). In a few other relevant contexts, LM was also investigated in density deconvolution in Comte et al.~(2008), and in the Laplace deconvolution in Benhaddou~(2018) where the unknown response function $f(\cdot)$ is non-periodic and defined on the entire positive real half-line. 

The objective of the paper is to extend the work of Benhaddou et al.~(2013) to the case when the noise is an anisotropic multi-parameter fractional Brownian sheet. Following a standard procedure, we derive minimax lower  bounds for the $L^p$-risk, with $1\leq p < \infty$, when $f(t , x)$ belongs to an anisotropic Besov ball and  the blurring function $g(t , x)$ is regular smooth. Taking advantage of the Riesz poly-potential operator, and inspired by the work of Wang~(1996, 1997), we apply the wavelet-vaguelette expansion (WVD) to de-correlate the anisotropic fBs in two directions. In addition, we take advantage of the flexibility of the Meyer wavelet basis in the Fourier domain to construct a wavelet hard-thresholding estimator that is adaptive and asymptotically quasi-optimal within a logarithmic factor of $\varepsilon$ over a large array of Besov balls. Furthermore, the proposed estimator attains convergence rates that depend on a delicate balance between the parameters of the Besov ball, the degree of ill-posedness of the convolution, and the parameters $\alpha_1$ and $\alpha_2$ associated with the fBs. Similar to the white noise case studied in Benhaddou et al.~(2013), the proposed approach is easily extended to recovering an $r$-dimensional function, $r> 2$, and the corresponding convergence rates turn out to be dimension-free. Finally, a simulation study is carried out to further confirm the results of our theoretical findings. 

The rest of the paper is organized as follows. Section 2 introduces some notation as well as the estimation algorithm. Section 3 describes the derivation of the lower bounds for the $L^p$-risk of estimators of $f$ as well as the upper bounds and establishes the asymptotic optimality of the estimator. Section 4 presents a limited simulations study to complement the theoretical results from Section 3. Section 5 extends the results in Sections 2 and 3 to the general $r$-dimensional case. Finally, Section 6 contains the proofs of the theoretical results.
 \section{Estimation Algorithm.}
 In what follows, denote the complex conjugate of $a$ by $\bar{a}$. 
 Let  $\tilde{Y}(m_1, m_2)$, $\tilde{B^{\alpha}}(m_1, m_2)$, $\tilde{g}(m_1, m_2)$ and $\tilde{f}(m_1, m_2)$ 
be the two-dimensional Fourier coefficients of functions $dY(t, x)$, $dB^{\alpha}(t, x)$, $g(t, x)$, and $f(t, x)$, respectively. \\
Consider a bandlimited wavelet basis $\psi_{j, k}(y)$ (e.g., Meyer-type) on $[0, 1]$, and form the product wavelet basis 
\be \label{bxb}
\Psi_{\omega}(t, x)= \psi_{j_1, k_1}(t)\psi_{j_2, k_2}(x), \ (t, x)\in U,
\ee
where $\omega \in \Omega$, and 
\be \label{omega}
\Omega=\left\{\omega=(j_1, k_1; j_2, k_2): j_i=m_0-1, \cdots, \infty; k_i=0, 1, \cdots, 2^{j_i}-1, i=1, 2 \right\}.
\ee 
Let $m_0$ be the lowest resolution level for the Meyer basis and denote the scaling function for the wavelet by $\psi_{m_0-1, k_i}$, $i=1, 2$. Since the functions \fr{bxb} form an orthonormal basis of the $L^2(U)$-space, the function $f(t, x)$ can be expanded over these bases with coefficients $\beta_{\omega}$ into a wavelet series 
\be \label{wavexp}
f(t, x)=\sum_{\omega \in \Omega}\beta_{\omega}\Psi_{\omega}(t, x),
\ee
where $\beta_{\omega}=\int^1_0\int^1_0f(t, x)\Psi_{\omega}(t, x)dtdx$. Denote the two-dimensional Fourier transform of $\Psi_{\omega}(t, x)$ by $\tilde{\Psi}_{\omega}(m_1, m_2)=\tilde{\psi}_{j_1, k_1}(m_1)\tilde{\psi}_{j_2, k_2}(m_2)$.  It is well-known (see, e.g, Johnstone et al~(2004), section 3.1) that under the Fourier domain and for any $j_i \geq m_0$, $i=1, 2$, one has  
\be \label{omeg}
W_{j_i}=\left \{ m_i: \tilde{\psi}_{j_i, k_i}(m_i)\neq 0 \right\} \subseteq 2\pi/3\left[ -2^{j_i+2}, -2^{j_i} \right] \cup \left[ 2^{j_i}, 2^{j_i+2} \right ].
\ee
Apply the two-dimensional Fourier transform to equation \fr{conveq} to obtain 
\be \label{fconv}
\tilde{Y}(m_1, m_2)=\overline{\tilde{f}(m_1, m_2)}\tilde{g}(m_1, m_2)+ \varepsilon^{\overline{\alpha}}\tilde{B}^{{\alpha}}(m_1, m_2).
\ee 
Now, applying Plancherel formula in both directions, the wavelet coefficients $\beta_{\omega}$ in \fr{wavexp} can be rewritten as 
\be \label{betinf}
\beta_{\omega}=\sum_{m_1\in W_{j_1}}\sum_{m_2 \in W_{j_2}}\overline{\tilde{\psi}_{j_1, k_1}(m_1)} \tilde{\psi}_{j_2, k_2}(m_2)\tilde{f}(m_1, m_2).
\ee 
Therefore, an unbiased estimator for $\beta_{\omega}$ is given by
\be \label{bethat}
\tilde{\beta}_{\omega}=\sum_{m_1\in W_{j_1}}\sum_{m_2 \in W_{j_2}}\overline{\tilde{\psi}_{j_1, k_1}(m_1)} \tilde{\psi}_{j_2, k_2}(m_2)\frac{\tilde{Y}(m_1, m_2)}{\tilde{g}(m_1, m_2)}.
\ee
Finally, consider the hard thresholding estimator for $f(t, x)$
\be \label{ef-hat}
\widehat{f}_{\varepsilon}(t, x)= \sum_{\omega \in \Omega(J_1, J_2)} \tilde{\beta}_{\omega}\II\left(|\tilde{\beta}_{\omega}|> \lambda^{\alpha}_{j;\varepsilon}\right)\Psi_{\omega}(t, x), 
\ee
where
\be \label{Omega}
\Omega(J_1, J_2)=\left\{\omega=(j_1, k_1; j_2, k_2): j_i=m_0-1, \cdots, J_i; k_i=0, 1, \cdots, 2^{j_i}-1, i=1, 2\right\},
\ee
and the quantities $J_1$, $J_2$ and $\lambda^{\alpha}_{j;\varepsilon}$ are to be specified. \\
\noindent
{\bf Assumption A.2.  }\label{A1} In the Fourier domain, the convolution kernel $g(t, x)$, for some positive constants $\nu$, and $C_1$ and $C_2$, independent of $m_1$ and $m_2$, is such that 
\be  \label{blur}
 C_1|m_1|^{-2\nu} \leq |\tilde{g}(m_1, m_2)|^2 \leq C_2 |m_1|^{-2\nu}.
\ee
Let us now evaluate the variance of \fr{bethat}.   
 \begin{lemma} \label{lem:Var}
Let $\tilde{ \beta}_{\omega}$ be defined in \fr{bethat}. Then, under the condition \fr{blur}, for   $1\leq p <\infty$, one has 
\be \label{var}
\EE \left|\tilde{ \beta}_{\omega}-\beta_{\omega}\right|^{2p} \leq K_1 \varepsilon^{2p\overline{\alpha}}2^{j_1p(2\nu+\alpha_1-1)+ j_2p(\alpha_2-1)}.
\ee
 In addition, 
\be \label{varp}
\EE \left|\tilde{ \beta}_{\omega}-\beta_{\omega}\right|^{p} \leq K_2 \varepsilon^{{p}\overline{\alpha}}2^{j_1\frac{p}{2}(2\nu+\alpha_1-1)+ j_2\frac{p}{2}(\alpha_2-1)}.
\ee
where $K_1$  and $K_2$ are some positive constants independent of $\varepsilon$.  
\end{lemma}
According to Lemma \ref{lem:Var}, choose the thresholds $\lambda^{\alpha}_{j;\varepsilon}$  of the form  
\be  \label{Thres}
\lambda^{\alpha}_{j;\varepsilon}=  \gamma  \varepsilon^{\overline{\alpha}} \sqrt{|\ln(\varepsilon)|}2^{{\frac{j_1}{2}(2\nu +\alpha_1-1)}}2^{{\frac{j_2}{2}(\alpha_2-1)}}, 
\ee 
 where $\gamma$ is some positive constant independent of $\varepsilon$. Furthermore, the finest resolution levels $J_1$ and $J_2$ are such that 
\be  \label{Lev:J}
2^{J_2}= \left[ \frac{\varepsilon^{2\overline{\alpha}}}{A^2}\right]^{-\frac{1}{\alpha_2}}, \ \ \ 2^{J_1}=  \left[ \frac{\varepsilon^{2\overline{\alpha}}}{A^2}\right]^{-\frac{1}{2\nu + \alpha_1}}.
\ee
\section{Convergence rates and asymptotic optimality.}
Denote 
\beqn  \label{eq10}
 s^{*}_i&=&s_i+1/2 - 1/\pi,\\
 s'_i&=&s_i + 1/2 -1/\pi',\\
 s''_i&=& s_i +1/p -1/p',
\eeqn
where $\pi'=\min\{2, \pi\}$ and $p'=\min\{p, \pi\}$. \\
\noindent
{\bf Assumption A.3.} The function $f(t, x)$ belongs to an anisotropic two-dimensional Besov space. In particular,  its wavelet coefficients $ \beta_{j_1, k_1, j_2, k_2}$ satisfy
\be  \label{eq11}
 B^{s_1, s_2}_{\pi, q}(A)=\left \{ f \in L^2(U): \left( \sum_{j_1, j_2} 2^{(j_1s_1^{*}+j_2s_2^*)q}\left (\sum_{k_1, k_2}| \beta_{j_1, k_1, j_2, k_2}|^{\pi}\right)^{q/{\pi}}\right )^{1/q} \leq A\right \} 
\ee
To construct minimax lower bounds for the $L^p$-risk, we define the $L^p$-risk over the set  $V$ as 
\be  \label{eq12}
 R(V)=\inf_{\tilde{f}} \sup _{f \in V}\EE \| \tilde{f}_{\varepsilon}-f\|^p_p,
\ee
where $\|g\|_p$ is the $L^p$-norm of a function $g$ and the infimum is taken over all possible estimators $\tilde{f}$ of $f$. 
\begin{theorem} \label{th:lowerbds}
Let $\min\{s_1, s_2\} \geq \max\{\frac{1}{\pi}, \frac{1}{2} \}$ with $1 \leq \pi,q \leq \infty$, and $A > 0$. Then, under conditions \fr{blur} and \fr{eq11}, for   $1\leq p <\infty$, as $\varepsilon \rightarrow 0$, 
 \be \label{lowerbds}
R(B^{s_1, s_2}_{\pi, q}(A)) \geq C A^p\left\{ \begin{array}{ll} 
 \left[\frac{ \varepsilon^{2\overline{\alpha}}}{A^2} \right]^{\frac{ps_2}{2s_2+\alpha_2}} , & \mbox{if}\ \ s_1 > \frac{s_2}{ \alpha_2}(2\nu + \alpha_1), \ \& \ \frac{s_2}{ \alpha_2}\geq \frac{p}{2}(\frac{1}{p'}-\frac{1}{p})\\
 \left[ \frac{ \varepsilon^{2\overline{\alpha}}}{A^{2}}\right]^{\frac{ps_1}{2s_1 +2\nu +\alpha_1}}   , & \mbox{if}\ \ \  \frac{p}{2}(2\nu + \alpha_1)(\frac{1}{\pi}-\frac{1}{p})\leq s_1\leq \frac{s_2}{\alpha_2}(2\nu + \alpha_1),\\
  \left[\frac{ \varepsilon^{2\overline{\alpha}}}{A^2} \right]^{\frac{p(s_1 + \frac{1}{p} -\frac{1}{\pi})}{2s^*_1+2\nu +\alpha_1-1}}, &  \mbox{if}\  s_1< (2\nu + \alpha_1)\frac{p}{2}(\frac{1}{\pi}-\frac{1}{p}), \ \& \ \frac{s_2}{ \alpha_2}\geq \frac{p}{2}(\frac{1}{p'}-\frac{1}{p}).
\end{array} \right.
\ee
  \end{theorem}
{\bf Proof of Theorem \ref{th:lowerbds}}. In order to prove the theorem, first we establish the equivalence between a process that involves the fractional Wiener sheet (eq. \fr{conveq}) and another that involves the standard Wiener sheet by applying an appropriate integral operator, taking advantage of relation \fr{fb-sb}. Then, we consider two cases, the case
when $f(t,x)$ is dense in both dimensions (the dense-dense case) and the case when $f(t,x)$ is dense in $x$ but sparse in $t$. Lemma $A.1$ of Bunea et al.~(2007) is then applied to find such lower bounds using conditions \fr{blur} and \fr{eq11}, combined with some useful properties of Meyer wavelet basis. To complete the proof, we choose the highest of the lower bounds. $\Box$ \\
The derivation of upper bounds of the $L^p$-risk relies on the following two lemmas.
\begin{lemma} \label{lemmabetap}
Under condition \fr{eq11}, one has
\be \label{betap}
\sum\limits_{k_1=0}^{2^{j_1} -1} \sum\limits_{k_2=0}^{2^{j_2} -1} \mid \beta_{j_1,k_1,j_2,k_2} \mid^p \leq A^p 2^{-p[(j_1s_1+j_2s_2)+ (\frac{1}{2} - \frac{1}{p'} )(j_1+j_2) ]} , \ \forall j_1, j_2\geq0.
\ee
\end{lemma}
\begin{lemma} \label{lem:Lar-D} 
Let $\tilde{ \beta}_{\omega}$ and $\lambda^{\alpha}_{j;\varepsilon}$ be defined by \fr{bethat} and \fr{Thres}, respectively. Define, for some positive constant $\eta$, the set 
\be  \label{Set:Thet}
\Theta_{\omega, \eta}=\left \{ \Theta: |\tilde{ \beta}_{\omega}-\beta_{\omega}| > \eta\lambda^{\alpha}_{j; \varepsilon} \right\}.
\ee
 Then, under condition \fr{blur} and as $\varepsilon \rightarrow 0$, one has 
\be  \label{Largdev}
\Pr \left( \Theta_{\omega, \eta}\right)= O\left(\frac{\left[\varepsilon^{2\overline{\alpha}}\right]^{\tau}}{{|\ln(\varepsilon)|^{\frac{1}{2}}}}\right),
\ee
where $\tau=\frac{\eta^2\gamma^2C^2_1}{8\overline{\alpha}}\left(\frac{3}{8\pi}\right)^{2\nu +\alpha_1-1+\alpha_2-1}$, $C_1$ and $\gamma$ appear in \fr{blur} and \fr{Thres}, respectively, and $\overline{\alpha}=\frac{\alpha_1+\alpha_2}{2}$.
\end{lemma}
\begin{theorem} \label{th:upperbds}
Let $\widehat{f}(. , .)$ be the wavelet estimator in \fr{ef-hat}, with $J_1$ and $J_2$ given by \fr{Lev:J}. Let  $\min\{s_1, s_2\} \geq \max\{\frac{1}{\pi}, \frac{1}{2} \}$ with $1 \leq \pi,q \leq \infty$, and let conditions  \fr{blur} and \fr{eq11} hold. If $\gamma$ in \fr{Thres} is large enough, then, for $1\leq p \leq\infty$ and as $\varepsilon \rightarrow 0$, one has
 \be \label{upperbds}
R( B^{s_1, s_2}_{\pi, q}(A)) \leq C A^p\left\{ \begin{array}{ll} 
  \left[\frac{ \varepsilon^{2\overline{\alpha}}}{A^2} \right]^{\frac{ps_2}{2s_2+\alpha_2}}\left[|\ln(\varepsilon)|\right]^{\xi_1+\frac{ps_2}{2s_2+\alpha_2}} , & \mbox{if}\  \ s_1 \geq \frac{s_2}{ \alpha_2}(2\nu + \alpha_1)\  \&\ \frac{s_2}{\alpha_2}\geq \frac{p}{2}(\frac{1}{p'}-\frac{1}{p}),\\
   \left[\frac{ \varepsilon^{2\overline{\alpha}}}{A^2} \right]^{\frac{ps_1}{2s_1 +2\nu +\alpha_1}}\left[|\ln(\varepsilon)|\right]^{\frac{ps_1}{2s_1 +2\nu +\alpha_1}} , & \mbox{if} \ \  \frac{p}{2}(2\nu + \alpha_1)(\frac{1}{\pi}-\frac{1}{p})< s_1< \frac{s_2}{\alpha_2}(2\nu + \alpha_1),\\
 \left[\frac{ \varepsilon^{2\overline{\alpha}}}{A^2} \right]^{\frac{p(s_1 + \frac{1}{p} -\frac{1}{\pi})}{2s^*_1+2\nu + \alpha_1-1}}\left[|\ln(\varepsilon)|\right]^{\xi_2+{\frac{p(s_1 + \frac{1}{p} -\frac{1}{\pi})}{2s^*_1+2\nu + \alpha_1-1}}}, &  \mbox{if}\ \  s_1\leq (2\nu + \alpha_1)\frac{p}{2}(\frac{1}{\pi}-\frac{1}{p})\ \& \ \frac{s_2}{\alpha_2} \geq\frac{p}{2}(\frac{1}{p'}-\frac{1}{p}),
\end{array} \right.
\ee
where $\xi_1$ and $\xi_2$ are defined as 
\beqn \label{d}
\xi_1 &=&\II \left(s_1= \frac{s_2}{ \alpha_2}(2\nu +\alpha_1) \right),\\
\xi_2 &=&\II  \left( s_1= (2\nu +\alpha_1)\frac{p}{2}\left(\frac{1}{\pi} -\frac{1}{p}\right) \right)+ \II\left(\frac{s_2}{\alpha_2}=\frac{p}{2}\left(\frac{1}{p'}-\frac{1}{p}\right)\right).
\eeqn
\end{theorem}
{\bf The proof of Theorem \ref{th:upperbds}.} The proof is very similar to that of Theorem 2 in Benhaddou.~(2017). $\Box$
\begin{remark}	
Notice that the finest resolution levels $J_1$ and $J_2$ in \fr{Lev:J} and the thresholds $\lambda^{\alpha}_{j;\varepsilon}$ in \fr{Thres}, are independent of the unknown parameters of the Besov ball  \fr{eq11} and therefore estimator \fr{ef-hat} is adaptive with respect to those parameters.
\end{remark}
\begin{remark}	
 Theorems \ref{th:lowerbds} and \ref{th:upperbds} imply that, for the $L^p$-risk, the estimator in \fr{ef-hat} is asymptotically quasi-optimal within a logarithmic factor of $\varepsilon$, over a wide range of anisotropic Besov balls $ B^{s_1, s_2}_{\pi, q}(A)$. 
\end{remark}
\begin{remark}	
The convergence rates depend on a delicate balance between the parameters of the Besov ball, smoothness of the convolution kernel $\nu$ and the parameters of the anisotropic fractional Brownian sheet, $\alpha_1$ and $\alpha_2$. These rates deteriorate as the values of $\alpha_1$ and $\alpha_2$ get closer and closer to zero.
\end{remark}
\begin{remark}	
To de-correlate the two-dimensional fractional Gaussian noise a wavelet-vaguelette expansion based on Meyer wavelet bases has been applied to the fractional Brownian sheet. The validity of a wavelet-based series representation for the fractional Brownian motion has been established in Wang~(1997), Meyer, Sellan and Taqqu~(1999) and Ayache and Taqqu~(2002). In particular, Ayache and Taqqu~(2002) established the optimality of wavelet-based series representation of the fractional Brownian sheets.
\end{remark}\begin{remark}	
For $\alpha_1=\alpha_2=1$, the rates of convergence match exactly those in Benhaddou~(2017), and, with $p=2$, those in Benhaddou et al.~(2013) in their bivariate white noise case.
\end{remark}
\begin{remark}	
If we hold the second dimension $x$ fixed, our rates are comparable to those in Wang~(1997) and Wishart~(2013) in their standard deconvolution with the one-parameter fractional Gaussian noise case.
\end{remark}
\section{Simulation Study}
In order to investigate the effect of the long-memory on the performance of our estimator, we carried out a limited simulation study.  We implemented an estimation algorithm for the model \fr{disconveq} using a modification of WaveD method of Raimondo and Stewart~(2007). We evaluated mean integrated square error (MISE) $\EE \| \hat{f}-f\|^2$ of the functional deconvolution estimator. 
\begin{enumerate}
\item We generated the data using equation \fr{disconveq} with convolution kernel $g(t,x) = 0.5 \exp(-|t| (1 + (x - 0.5)^2))$, and $f(t,x)=f(t)f(x)$. In particular, we chose $f(t)$ to be a LIDAR or Doppler signal over the grid $t_i=\frac{i}{N}$ with $i = 1,2, . . . ,N$, and we chose $f(x)=\exp(-|x-0.5|x^3)$ over the grid $x_i=\frac{i}{N}$ with $i = 1,2, . . . ,N$, and $N=2^{10}=1024$. The LIDAR and Doppler signals were generated from the {\fontfamily{qcr}\selectfont \small waved} package. All test functions were scaled to have a unit norm. Bear in mind that, though $f(t,x)$ is a product of two univariate functions, the method does not ``recognize'' this and, therefore, cannot take advantage of this information. Also, notice that with our choice of $g(t, x)$, the degree of ill-posedness (DIP) of the convolution is DIP=0.5.
\item We simulated the LM error in the direction of both $t$ and $x$ using the {\fontfamily{qcr}\selectfont \small fracdiff} package of {\fontfamily{qcr}\selectfont \small R} available from {\fontfamily{qcr}\selectfont \small CRAN}. In particular, the {\fontfamily{qcr}\selectfont \small fracdiff.sim} command was used to simulate two one-dimensional fractionally differenced ARIMA$(0,d,0)$ (fARIMA) sequences which behave similar to a fractional Gaussian noise, and then multiply them together to generate a two-dimensional error structure. For the LM parameters $\alpha=(\alpha_1, \alpha_2)$, we used different combinations of the levels $0.8$, $0.6$, $0.4$ and $0.2$. To create a dependence structure that is anisotropic, we only used combinations such that $\alpha_1\neq \alpha_2$. Note that the fractional differencing parameter $d$ is obtained from $\alpha$ by the relation $d=\frac{1-\alpha}{2}$ . 
\item The choice of $\sigma$ in \fr{disconveq} was determined by the blurred signal-to-noise ratio (SNR), where
$
\text{SNR}=10\log_{10}\left( \frac{ \| q(t, x)  \|^2   }{\sigma^2}\right).
$
We considered three choices, SNR=10dB (high noise), 20dB (medium noise) and 30dB (low noise).
\item To compute the estimator, we used the function {\fontfamily{qcr}\selectfont \small FWaveD} and {\fontfamily{qcr}\selectfont \small IWaveD} from the {\fontfamily{qcr}\selectfont \small waved} package. First, we applied the Meyer wavelet transform to the data in $t$-direction by {\fontfamily{qcr}\selectfont \small FWaveD}, with no thresholding by setting {\fontfamily{qcr}\selectfont \small thr} to be a zero vector. Then, we applied the second thresholded Meyer wavelet transform in $x$-direction by {\fontfamily{qcr}\selectfont \small FWaveD}, however, the level dependent threshold were modified according to \fr{Thres}. For the tuning parameter $\gamma$, we tried different values ranging from $\gamma=\sqrt{2}$ to $\gamma=\sqrt{6}$, but the performance of the estimator reached its best at the latter value. Therefore, our default choice of $\gamma$ was $\sqrt{6}$. The finest resolution level $J_1$ was estimated using the default method in the {\fontfamily{qcr}\selectfont \small waved} package, together with the choice $J_2=5$. 
\end{enumerate}
Table 1 reports the averages of the errors over 2000 simulation runs using a sample size of $n=N^2= (2^{10})^2$. Figure \fr{fig1} illustrates the function $f(t.x)$ and its estimator $\widehat{f}(t, x)$ and shows a relatively good precision that deteriorates as the anisotropic pair of LM parameters $(\alpha_1, \alpha_2)$ decreases from 0.8 to 0.2. Notice here that according to Table 1, the influence of LM  is a little more pronounced in the direction of the dimension $x$ than it is in the direction of $t$. Table 1 complements Figure \fr{fig1} and confirms our theoretical results in the sense that as the level of LM increases ($\alpha$ gets smaller), the mean integrated squared error increases (the performance deteriorates).
\begin{figure}[htbp]
\caption{The function $f(t, x)$ (top) and its estimate $\widehat{f}(t, x)$ (bottom) at different combinations of $(\alpha_1, \alpha_2)$.\label{fig1}}
\scalebox{0.26}{
\begin{tabular}{ccc}
\includegraphics{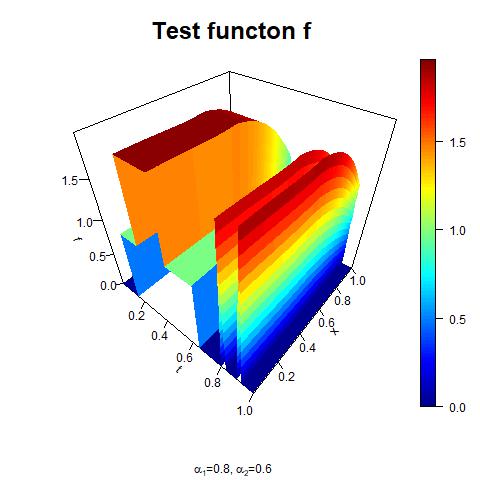}        &  \includegraphics{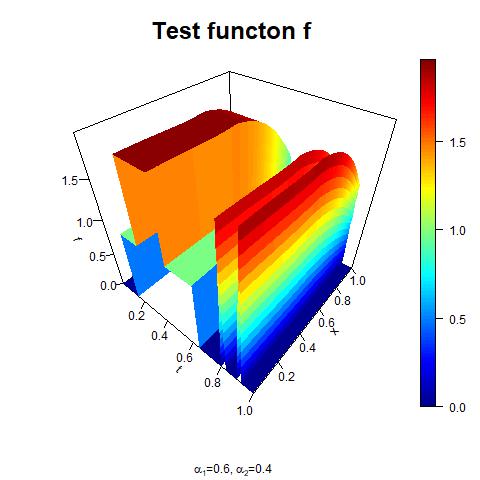}     &\includegraphics{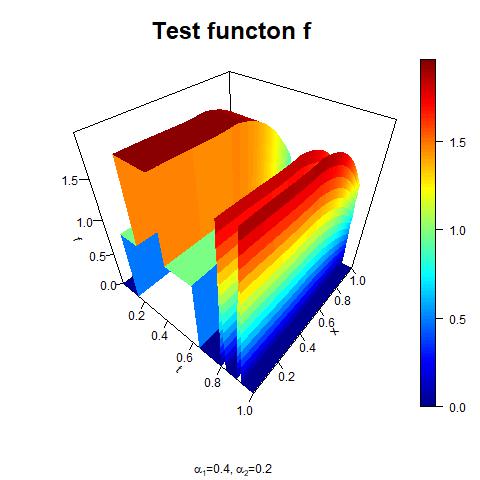}\\
\includegraphics{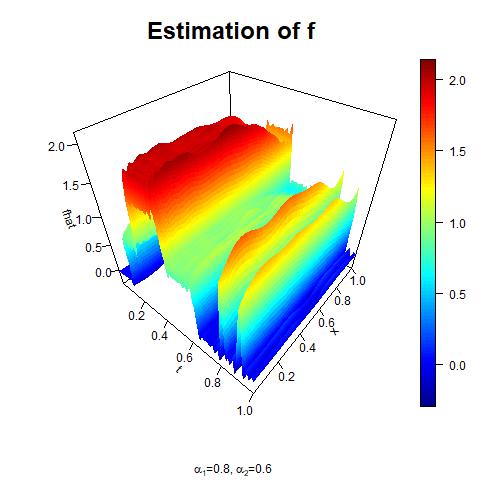}  &  \includegraphics{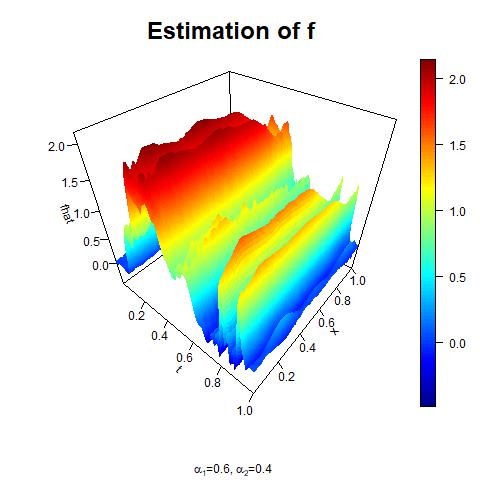}&\includegraphics{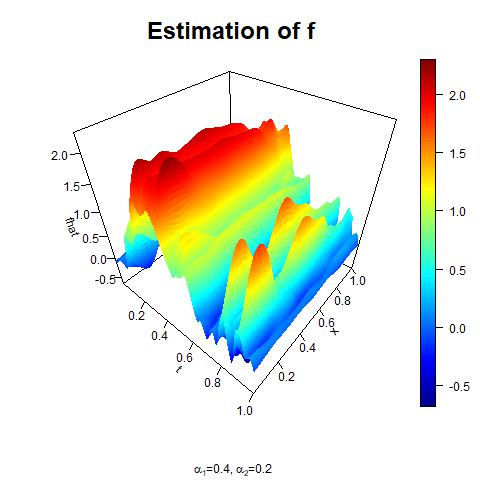}
\end{tabular}}
\end{figure}
\begin{table}[htbp]\centering
\caption{}
Finest resolution level $J_1$ are listed in the parentheses.
\scalebox{0.85}{
\begin{tabular}{|c|c|cccccc|}
\hline
DIP=0.5&&($ \alpha_1 $,$ \alpha_2 $)& &&&&  \\
signal&$\sigma$              &(0.8,0.6)    &(0.8,0.4)     &(0.8,0.2)    & (0.6,0.4)   & (0.6,0.2)       &(0.4,0.2)       \\

\hline
Lidar (10dB)&  0.20         &0.0779 (3) &0.0948 (3) &0.1296 (3) &0.1134 (3)  &0.2097 (3) &0.3232 (3)  \\
Lidar (20dB)&  0.06         &0.0434 (4) &0.0498 (4) & 0.0641 (4) &0.0572 (4)  &0.0796 (4)&0.0948 (4)           \\
Lidar (30dB)&  0.02         &0.0132 (5) &0.0183 (5) &0.0292 (5)&0.0219 (5)    &0.0358 (5)&0.0444 (5)   \\
\hline 
Doppler (10dB)&  0.15    &0.2613 (3)&0.2724 (3)&0.3160 (3)&0.2935 (3)  &0.3828 (3)&0.4823 (3)  \\
Doppler (20dB)&  0.05    &0.1254 (4)&0.1298 (4)&0.1437 (4)&0.1371 (4)  &0.1584 (4)&0.1867 (4)   \\
Doppler (30dB)&  0.015  &0.0516 (5)&0.0530 (5) &0.0568 (5)&0.0544 (5)&0.0610 (5)&0.0667 (5)  \\
\hline
\end{tabular}}
\end {table}
\begin{table}[htbp]\centering
\scalebox{0.85}{
\begin{tabular}{|c|c|cccccc|}
\hline
DIP=0.5&&($ \alpha_1 $,$ \alpha_2 $)& &&&&  \\ 
signal&$\sigma$               &(0.6,0.8)  &(0.4,0.8)       &(0.2,0.8)          &(0.4,0.6)   &(0.2,0.6)       &(0.2,0.4)   \\

\hline
Lidar (10dB)&  0.20          &0.0760 (3) &0.0851 (3) &0.0956 (3)     &0.1038 (3)  &0.1293 (3)&0.2251 (3)  \\
Lidar (20dB)&  0.06          &0.0423 (4) &0.0453 (4) &0.0486 (4)    &0.0536 (4)   &0.0612 (4)&0.0791 (4)            \\
Lidar (30dB)&  0.02          &0.0120 (5) &0.0132 (5) &0.0149 (5)     &0.0180 (5)  &0.0210 (5)&0.0331 (5) \\
\hline 
Doppler (10dB)&  0.15    &0.2601 (3)&0.2640 (3)&0.2749 (3)       &0.2836 (3)  &0.3090 (3)&0.3922 (3) \\
Doppler (20dB)&  0.05    &0.1247 (4)&0.1264 (4)&0.1292 (4)       &0.1331 (4)  &0.1376 (4)&0.1578 (4)   \\
Doppler (30dB)&  0.015  &0.0513 (5)&0.0510 (5)&0.0522 (5)      &0.0527 (5)  &0.0541 (5)&0.0591 (5)     \\
\hline
\end{tabular}}
\end {table}

\section{The $r$-dimensional case: estimation, convergence rates and asymptotic optimality.}
Consider the $r$-dimensional case of model \fr{conveq}
\be \label{conveq-r}
dY(t, {\bf x})=q(t, {\bf x})dtd{\bf x} + \varepsilon^{\overline{\alpha}}dB^{\bf{\alpha}}(t, {\bf x}), \ q(t, {\bf x})=\int^1_0f(s, {\bf x})g(t-s, {\bf x})ds.
\ee
Here, $U=[0, 1]^r$, ${\bf x}=(x_1, x_2, \cdots, x_{r-1})$, $g(t , {\bf x})$ is the convolution kernel, $B^{{\alpha}}(t, {\bf x})$ is an anisotropic $r$-dimensional fBs, $\overline{\alpha}={\sum^r_{k=1}\alpha_k}/{r}$, and $\alpha_k=2-2H_k\in(0, 1]$, $k=1,2, \cdots, r$, is the parameter of the long-memory in the direction of the $k^{th}$ variable. An $r$-parameter fBs is defined as a zero-mean Gaussian process $\{ B^{{\alpha}}({\bf t}): {\bf t}\in[0,1]^r\}$ whose covariance function is of the form
\be \label{cov}
\Cov(B^{{\alpha}}({\bf t_1}), B^{{\alpha}}({\bf t_2})) =c\Pi^r_{k=1}\left[|t_{k1}|^{2H_k}+|t_{k2}|^{2H_k}-|t_{k1}-t_{k2}|^{2H_k}\right],
\ee 
for some $c\in \mathbb{R}$, $(t_{k1}, t_{k2})\in [0, 1]^2$, $k=1, 2, \cdots, r$, where $H_k=1-\alpha_k/2$ is the Hurst parameter in the direction of the $k^{th}$ dimension.\\
Consider a bandlimited wavelet basis $\psi_{j, k}(y)$ (e.g., Meyer-type) on $[0, 1]$, and form the product wavelet basis 
\be \label{bxb-r}
\Psi_{\omega_r}(t, {\bf x})= \psi_{j_1, k_1}(t)\Pi^r_{i=2}\psi_{j_i, k_i}(x_{i-1}), \ (t, {\bf x})\in U,
\ee
where $\omega_r \in \Omega_r$, and 
\be \label{omegar}
\Omega_r=\left\{\omega_r=(j_1, k_1; j_2, k_2; \cdots; j_r, k_r): j_i=m_0-1, \cdots, \infty; k_i=0, 1, \cdots, 2^{j_i}-1, i=1, 2, \cdots, r \right\}.
\ee 
Then, the function $f(t, {\bf x})$ can be expanded over these bases with coefficients $\beta_{\omega_r}$ into a wavelet series 
\be \label{wavexpr}
f(t, {\bf x})=\sum_{\omega_r \in \Omega_r}\beta_{\omega_r}\Psi_{\omega_r}(t, {\bf x}),
\ee
where $\beta_{\omega_r}=\int_Uf(t, {\bf x})\Psi_{\omega_r}(t, {\bf x})dtd{\bf x}$. Apply the $r$-dimensional Fourier transform to equation \fr{conveq-r} to obtain 
\be \label{fconv}
\tilde{Y}({\bf m_r})=\tilde{f}({\bf m_r})\tilde{g}({\bf m_r})+ \varepsilon^{\overline{\alpha}}\tilde{B}^{{\alpha}}({\bf m_r}),
\ee 
where ${\bf m_r}=(m_1, m_2, \cdots, m_r)$. Then, applying Plancherel formula in all directions, the wavelet coefficients $\beta_{\omega_r}$ in \fr{wavexpr} can be rewritten as 
\be \label{betinfr}
\beta_{\omega_r}=\sum_{m_1\in W_{j_1}}\sum_{m_2 \in W_{j_2}}\cdots \sum_{m_r \in W_{j_r}} \tilde{\Psi}_{\omega_r}({\bf m_r})\tilde{f}({\bf m_r}).
\ee 
Therefore, an unbiased estimator for $\beta_{\omega_r}$ is given by
\be \label{bethat-r}
\tilde{\beta}_{\omega_r}=\sum_{m_1\in W_{j_1}}\sum_{m_2 \in W_{j_2}}\cdots \sum_{m_r \in W_{j_r}} \tilde{\Psi}_{\omega_r}({\bf m_r})\frac{\tilde{Y}({\bf m_r})}{\tilde{g}({\bf m_r})}.
\ee
Finally, consider the hard thresholding estimator for $f(t, {\bf x})$
\be \label{ef-hatr}
\widehat{f}_{\varepsilon}(t, {\bf x})= \sum_{\omega_r \in \Omega_r({\bf J})} \tilde{\beta}_{\omega_r}\II\left(|\tilde{\beta}_{\omega_r}|> \lambda^{\alpha}_{{\bf j};\varepsilon}\right)\Psi_{\omega_r}(t, {\bf x}), 
\ee
where
\be \label{Omega}
\Omega_r({\bf J})=\left\{\omega_r=(j_1, k_1; j_2, k_2; \cdots; j_r, k_r): j_i=m_0-1, \cdots, J_i; k_i=0, 1, \cdots, 2^{j_i}-1, i=1, 2, \cdots, r\right\},
\ee
with ${\bf j}=(j_1, j_2, \cdots, j_r)^T$, and the quantities $J_1$, $J_2$, $\cdots$, $J_r$ and $\lambda^{\alpha}_{{\bf j};\varepsilon}$ are to be determined. \\
\noindent
{\bf Assumption A.4.  }\label{A4} In the Fourier domain,  the convolution kernel $g(t, {\bf x})$, for some positive constants $\nu$, and $C_{r1}$ and $C_{r2}$, independent of ${\bf m_r}=(m_1, m_2, \cdots, m_r)$, is such that 
\be  \label{blur-r}
 C_{r1}|m_1|^{-2\nu} \leq |\tilde{g}({\bf m_r})|^2 \leq C_{r2} |m_1|^{-2\nu}.
\ee
 \begin{lemma} \label{lem:Var_r}
Let $\tilde{ \beta}_{\omega_r}$ be defined in \fr{bethat-r}. Then, under the condition \fr{blur-r}, for   $1\leq p <\infty$, one has 
\be \label{var_r}
\EE \left|\tilde{ \beta}_{\omega_r}-\beta_{\omega_r}\right|^{2p} \leq K_3 \varepsilon^{p2\overline{\alpha}}2^{j_1p(2\nu+\alpha_1-1)+ j_2p(\alpha_2-1)+\cdots+j_{r}p(\alpha_{r}-1)}.
\ee
 In addition, 
\be \label{varp_r}
\EE \left|\tilde{ \beta}_{\omega_r}-\beta_{\omega_r}\right|^{p} \leq K_4 \varepsilon^{{p}\overline{\alpha}}2^{j_1\frac{p}{2}(2\nu+\alpha_1-1)+ j_2\frac{p}{2}(\alpha_2-1)+\cdots+ j_{r}\frac{p}{2}(\alpha_{r}-1)   }.
\ee
where $K_3$  and $K_4$ are some positive constants independent of $\varepsilon$.  
\end{lemma}
Consequently, choose the thresholds $\lambda^{\alpha}_{{\bf j};\varepsilon}$ of the from 
\be  \label{Thres-r}
\lambda^{\alpha}_{{\bf j};\varepsilon}=  \rho  \varepsilon^{\overline{\alpha}} \sqrt{|\ln(\varepsilon)|}2^{j_1\nu}\Pi^r_{i=1}2^{{\frac{j_i}{2}(\alpha_i-1)}}, 
\ee 
 where $\rho$ is some positive constant that is independent of $\varepsilon$, $\overline{\alpha}=\frac{1}{r}{\sum^r_{k=1}\alpha_k}$ and $J_1$ and $J_i$, $i=2, 3, \cdots, r$,  are such that 
\be  \label{Lev:Jr}
2^{J_1}=  \left[ \frac{\varepsilon^{2\overline{\alpha}}}{A^2}\right]^{-\frac{1}{2\nu + \alpha_1}}, \ \ \ 2^{J_i}= \left[ \frac{\varepsilon^{2\overline{\alpha}}}{A^2}\right]^{-\frac{1}{\alpha_i}}.
\ee
\begin{lemma} \label{lem:Lar-D-r} 
Let $\tilde{ \beta}_{\omega_r}$ and $\lambda^{\alpha}_{{\bf j};\varepsilon}$ be defined by \fr{bethat-r} and \fr{Thres-r}, respectively. Define, for some positive constant $\eta$, the set 
\be  \label{Set:Thet}
\Theta^r_{\omega_r, \eta}=\left \{ \Theta: |\tilde{ \beta}_{\omega_r}-\beta_{\omega_r}| > \eta\lambda^{\alpha}_{{\bf j}; \varepsilon} \right\}.
\ee
 Then, under condition \fr{blur-r}, as $\varepsilon \rightarrow 0$, one has 
\be  \label{Largdev}
\Pr \left( \Theta^r_{\omega_r, \eta}\right)= O\left(\frac{\left[\varepsilon^{{2\overline{\alpha}}}\right]^{\tau_r}}{{|\ln(\varepsilon)|^{\frac{1}{2}}}}\right),
\ee
where $\tau_r=\frac{\eta^2\rho^2C^2_{r1}}{8\overline{\alpha}}\left(\frac{3}{8\pi}\right)^{2\nu +\sum^r _{i=1}(\alpha_i-1)}$, and $C_{r1}$ and $\rho$ appear in \fr{blur-r} and \fr{Thres-r}, respectively.
\end{lemma}
\noindent
{\bf Assumption A.5.  } The $r$-dimensional function $f(t, {\bf x})$ belongs to an anisotropic Besov space. In particular,  its wavelet coefficients $ \beta_{j_1, k_1; j_2, k_2; \cdots ; j_r, k_r}$ satisfy
\be  \label{eq11-r}
 B^{s_1, s_2, \cdots, s_r}_{\pi, q}(A)=\left \{ f \in L^2(U): \left( \sum_{j_1, j_2, \cdots, j_r} 2^{{\bf j}^T{\bf s^*}q}\left (\sum_{k_1, k_2, \cdots, k_r}| \beta_{j_1, k_1; j_2, k_2; \cdots; j_r, k_r}|^{\pi}\right)^{q/{\pi}}\right )^{1/q} \leq A\right \} 
\ee
where ${\bf j}=(j_1, j_2, \cdots, j_r)^T$, ${\bf s^*}=(s^*_1, s^*_2, \cdots, s^*_r)^T$ and $s^*_i$ are defined in \fr{eq10}.\\
Define the quantity 
\be \label{s-op}
\frac{s_{2o}}{\alpha_{2o}}=\min\left\{ \frac{s_2}{\alpha_2}, \frac{s_3}{\alpha_3}, \cdots, \frac{s_r}{\alpha_r}\right\},
\ee
and denote the dimension of $f$ that corresponds to \fr{s-op} by $i_o$.
\begin{theorem} \label{th:lowerbds-r}
Let $\min\{s_1, s_2\} \geq \max\{\frac{1}{\pi}, \frac{1}{2} \}$ with $1 \leq \pi,q \leq \infty$, and $A > 0$. Then, under conditions \fr{blur-r} and \fr{eq11-r}, for   $1\leq p <\infty$, as $\varepsilon \rightarrow 0$, 
 \be \label{lowerbds}
R(B^{s_1, s_2, \cdots, s_r}_{\pi, q}(A)) \geq C A^p\left\{ \begin{array}{ll} 
 \left[\frac{ \varepsilon^{2\overline{\alpha}}}{A^2} \right]^{\frac{ps_{2o}}{2s_{2o}+\alpha_{2o}}} , & \mbox{if}\ \  s_1 > \frac{s_{2o}}{\alpha_{2o}}(2\nu + \alpha_1), \ \frac{s_{2o}}{\alpha_{2o}}\geq \frac{p}{2}(\frac{1}{p'}-\frac{1}{p})\\
 \left[ \frac{ \varepsilon^{2\overline{\alpha}}}{A^{2}}\right]^{\frac{ps_1}{2s_1 +2\nu +\alpha_1}}   , & \mbox{if}\ \ \  \frac{p}{2}(2\nu + \alpha_1)(\frac{1}{\pi}-\frac{1}{p})\leq s_1\leq \frac{s_{2o}}{\alpha_{2o}}(2\nu + \alpha_1),\\
  \left[\frac{ \varepsilon^{2\overline{\alpha}}}{A^2} \right]^{\frac{p(s_1 + \frac{1}{p} -\frac{1}{\pi})}{2s^*_1+2\nu +\alpha_1-1}}, &  \mbox{if}\  s_1< (2\nu + \alpha_1)\frac{p}{2}(\frac{1}{\pi}-\frac{1}{p}), \ \frac{s_{2o}}{\alpha_{2o}}\geq \frac{p}{2}(\frac{1}{p'}-\frac{1}{p}).
\end{array} \right.
\ee
 \end{theorem}
 \begin{theorem} \label{th:upperbds-r}
Let $\widehat{f}(. , .)$ be the wavelet estimator in \fr{ef-hatr}, with $J_1$ and $J_i$ with $i=2, \cdots, r$ given by \fr{Lev:Jr}. Let  $\min\{s_1, s_2, \cdots, s_r\} \geq \max\{\frac{1}{\pi}, \frac{1}{2} \}$ with $1 \leq \pi,q \leq \infty$, and let conditions  \fr{blur-r} and \fr{eq11-r} hold. If $\rho$ in \fr{Thres-r} is large enough, then, for $1\leq p \leq\infty$, as $\varepsilon \rightarrow 0$, 
 \be \label{upperbds}
R( B^{s_1, ., s_r}_{\pi, q}(A)) \leq C A^p\left\{ \begin{array}{ll} 
  \left[\frac{ \varepsilon^{2\overline{\alpha}}}{A^2} \right]^{\frac{ps_{2o}}{2s_{2o}+\alpha_{2o}}}\left[|\ln(\varepsilon)|\right]^{\xi_1+\frac{ps_{2o}}{2s_2+\alpha_{2o}}} , & \mbox{if}\  \ s_1 \geq \frac{s_{2o}}{ \alpha_{2o}}(2\nu + \alpha_1), \ \frac{s_{2o}}{\alpha_{2o}}\geq \frac{p}{2}(\frac{1}{p'}-\frac{1}{p}),\\
   \left[\frac{ \varepsilon^{2\overline{\alpha}}}{A^2} \right]^{\frac{ps_1}{2s_1 +2\nu +\alpha_1}}\left[|\ln(\varepsilon)|\right]^{ \frac{ps_1}{2s_1 +2\nu +\alpha_1}} , & \mbox{if} \ \  \frac{p}{2}(2\nu + \alpha_1)(\frac{1}{\pi}-\frac{1}{p})< s_1< \frac{s_{2o}}{\alpha_{2o}}(2\nu + \alpha_1),\\
 \left[\frac{ \varepsilon^{2\overline{\alpha}}}{A^2} \right]^{\frac{p(s_1 + \frac{1}{p} -\frac{1}{\pi})}{2s^*_1+2\nu + \alpha_1-1}}\left[|\ln(\varepsilon)|\right]^{\xi_2+{\frac{p(s_1 + \frac{1}{p} -\frac{1}{\pi})}{2s^*_1+2\nu + \alpha_1-1}}}, &  \mbox{if}\ \  s_1\leq (2\nu + \alpha_1)\frac{p}{2}(\frac{1}{\pi}-\frac{1}{p}), \ \frac{s_{2o}}{\alpha_{2o}}\geq \frac{p}{2}(\frac{1}{p'}-\frac{1}{p}),
\end{array} \right.
\ee
where $\xi_1$ and $\xi_2$ are defined as 
\beqn \label{d}
\xi_1 &=&\II \left(s_1= \frac{s_{2o}}{ \alpha_{2o}}(2\nu +\alpha_1) \right) +\sum^r_{i=2, i\neq i_o}\II\left(\frac{s_{2o}}{\alpha_{2o}}=\frac{s_i}{\alpha_i}\right),\\
\xi_2 &=&\II  \left( s_1= (2\nu +\alpha_1)\frac{p}{2}\left(\frac{1}{\pi} -\frac{1}{p}\right) \right)+\sum^r_{i=2, i\neq i_o}\II\left(\frac{p}{2}\left(\frac{1}{p'}-\frac{1}{p}\right)=\frac{s_i}{\alpha_i}\right).
\eeqn
\end{theorem}
{\bf Proof of Theorem \ref{th:upperbds-r}.} The proof is very similar to that of Theorem 5 in Benhaddou et al.~(2013) so we skip it. $\Box$
 \begin{remark}	
 Notice that the convergence rates in Theorems \ref{th:lowerbds-r} and \ref{th:upperbds-r} depend on a delicate balance between $s_1$, $\pi$, $\nu$, $\alpha_1$, $\min\left\{ \frac{s_i}{\alpha_i}, i=2, 3, \cdots, r\right\}$ and the average of the parameters of the anisotropic fractional Brownian sheet, $\alpha_i$, $i=1, 2, \cdots, r$. Such rates are comparable to those in Benhaddou et al.~(2013) in their white noise case when $\alpha_1=\alpha_2=\cdots= \alpha_r=1$. In addition, the rates are independent of the dimension $r$ as expected. 
\end{remark}

\section{Proofs. }
In the proofs of Theorems \ref{th:lowerbds} and \ref{th:upperbds} we will use the properties described in Petsa and Sapatinas~(2009) associated with Meyer wavelet basis. Namely, the properties of concentration, unconditionality and Temlyakov.
\subsection{Proof of the lower bounds.}
In order to prove Theorem \ref{th:lowerbds}, we use Lemma A.1 
of \cite{bun}.
\begin{lemma} \label{lem:Bunea} 
Let $\Te$ be a set of functions of cardinality $\card(\Te)\geq 2$ such that\\
(i) $\|f-g\|_p^p \geq 4\delta^p, \ for\  f, g \in \Te, \ f \neq g, $\\
(ii) the Kullback divergences $K(P_f, P_g)$ between the measures $P_f$ and $P_g$ 
satisfy the inequality $K(P_f, P_g) \leq \log(\card(\Te))/16,\ for\ f,\ g \in \Te$.\\
Then, for some absolute positive constant $C_1$, one has 
$$
 \inf_{f_n}\sup_{f\in \Te} \EE_f\|f_n-f\|_p^p \geq C_1 \delta^p,
$$
where $\inf_{f_n}$ denotes the infimum over all estimators.
\end{lemma}
{\bf Proof of Theorem \ref{th:lowerbds}}. Let $P_f$ be the probability law of the process $\{q(t, x) + \varepsilon^{\overline{\alpha}} dB^{\alpha}(t, x), (t, x)\in U\}$ when $f$ is true. Let us first transform this process into a process that involves the standard Brownian sheet, taking into account relation \fr{fb-sb}. This can be accomplished using the following functions
\be \label{gammafu}
\Gamma_{\alpha}(t, u)= \frac{C_H}{\Gamma(H+1/2) \Gamma(3/2-H)}\left[(t-u)_{+}\right]^{1/2-H},\ for\ t\in [0, 1],
\ee
where $C_H$ is some explicit constant, $\Gamma(\cdot)$ is the Gamma function. Consider $0< u < t < 1$, then \fr{gammafu} reduces to 
\be \label{gammafu2}
\Gamma_{\alpha}(t, u)= C(t-u)^{1/2-H}.
\ee
Now define the operator $K_{\Gamma}$ such that
\be \label{kappagamma}
K_{\Gamma}f(t, x)= \int_R\int_R \Gamma_{\alpha_1}(t, u) \Gamma_{\alpha_2}(x, v)f(u, v)dudv.
\ee
Taking into account relation \fr{fb-sb}, it is well known that operator \fr{kappagamma} converts an anisotropic two-dimensional fractional Brownian sheet with Hurst parameters $(H_1, H_2)=(1-{\alpha_1}/{2},1-\alpha_2/2)$ into a standard Brownian sheet. Let us then apply operator \fr{kappagamma} to equation \fr{conveq} to obtain, 
\be \label{haze2}
\widehat{Y}(t, x)=K_{\Gamma}dY(t, x)=K_{\Gamma}q(t, x) + \varepsilon^{\overline{\alpha}} B(t, x), 
\ee
where $B(t, x)$ are standard Brownian sheet. Finally, differentiating both sides of \fr{haze2} with respect to $t$ and $x$, yields
\be \label{haze3}
d\widehat{Y}(t, x)=K'_{\Gamma}q(t, x)dtdx +  \varepsilon^{\overline{\alpha}} dB(t, x),
\ee
where $K'_{\Gamma}q(t, x)$ is 
\be\label{derkap}
K'_{\Gamma}q(t, x)=\frac{d^2}{dtdx}K_{\Gamma}q(t, x)=c\int^x_0(x-v)^{-1/2-H_2}\int^t_0f(s, v)\int^{t-s}_0(t-s-r)^{-1/2-H_1}g(r, v)drdsdv.
\ee
Notice that in the Fourier domain, by the convolution theorem and condition \fr{blur}, \fr{derkap} has the form
\be \label{fderkap}
\tilde{K'_{\Gamma}}q(m_1, m_2)\asymp \tilde{f}(m_1, m_2)|m_1|^{-\nu}|m_1|^{(1-\alpha_{1})/2}|m_2|^{(1-\alpha_{2})/2}.
\ee
\underline{\bf  The dense-dense case.}
Let $\eta$ be the matrix with components $\eta_{k_1k_2}=\{0, 1\}$, $k_i=0, 1, \cdots, 2^{j_i}-1$, $i=1, 2$, and denote the set of all possible values $\eta$ by $H$. Define the functions
\be 
f_{j_1j_2}(t, x)=\gamma_{j_1j_2}\sum^{2^{j_1}-1}_{k_1=0}\sum^{2^{j_2}-1}_{k_2=0}\eta_{k_1k_2}\psi_{j_1k_1}(t)\psi_{j_2k_2}(x).\label{tesf1}
\ee
Note that the matrix $\eta$ has $N=2^{j_1+j_2}$ components and therefore, $\card(H)=2^N$. To guarantee that functions \fr{tesf1} satisfy \fr{eq11}, choose $\gamma_{j_1j_2}=A 2^{-j_1(s_1+1/2)}2^{-j_2(s_2+1/2)}$. Now, take another function of the form \fr{tesf1}, $\tilde{f}_{j_1j_2}$, with $\tilde{\eta}_{k_1k_2} \in H$ instead of $\eta_{k_1k_2}$ and compute the $L^p$-norm of the difference between $f_{j_1j_2}$ and $\tilde{f}_{j_1j_2}$, with the application of the Varshamov-Gilbert Lemma (\cite{tsybakov}, p 104), to obtain   
\be  \label{Norm}
\|\tilde{f}_{j_1j_2}-f_{j_1j_2}\|^p_p \geq \gamma^p_{j_1j_2} 2^{j_1p/2}2^{j_2p/2}/8.
\ee
In addition, define $\widehat{P}_f$, the probability law of the process $\left\{K'_{\Gamma}q(t, x)+ \varepsilon^{\overline{\alpha}}dW(t, x), (t, x)\in U\right\}$ when $f$ is true, where $W(t, x)$ is a two-parameter standard Wiener sheet. Note that the probability measures $\widehat{P}_f$ and ${P}_f$ are stochastically equivalent (see, e.g., Huang et al.~(2006), Theorem 3.1). Hence, by the multi-parameter Girsanov formula (see, e.g., Dozzi~(1989), p.89), the Kullback divergence takes the form  
\be  \label{KLbk}
K({P}_f, {P}_{{\tilde{f}}})=K(\widehat{P}_f, \widehat{P}_{{\tilde{f}}})= \EE \left[ \log\left( \widehat{P}_f/ \widehat{P}_{\tilde{f}}\right)\right]_{\widehat{P}_f}= \left(2\varepsilon^{2\overline{\alpha}} \right)^{-1}\left|\left|{K'_{\Gamma}g*(\tilde{f}-f)(t, x)}{}\right|\right|^2.
\ee
Since $\left|\tilde{\eta}_{k_1k_2}-\eta_{k_1k_2}\right| \leq 1$, plugging $f_{j_1j_2}$ and $\tilde{f}_{j_1j_2}$ into \fr{KLbk}, applying Plancherel's formula, $|\psi_{j, k, m}| \leq 2^{-j/2}$ and \fr{fderkap}, it can be shown that 
\beqn  \label{KLbk:upperbd}
K(\widehat{P}_f, \widehat{P}_{\tilde{f}}) &\leq& C \left(2\varepsilon^{2\overline{\alpha}} \right)^{-1}{\gamma^2_{j_1j_2}}2^{-j_1(2\nu +\alpha_1-1)}2^{-j_2(\alpha_2-1)}2^{j_1 + j_2}
.\eeqn
Part $(ii)$ of Lemma \ref{lem:Bunea} gives the constraint
\be
2^{-j_1(2s_1 + 2\nu +\alpha_1)}2^{-j_2(2s_2 +\alpha_2)} \leq \frac{C \varepsilon^{\alpha_1+\alpha_2}}{A^2}.
\ee
Now, define 
\be
\tau_{\varepsilon}=\log_2\left(\frac{CA^2}{\varepsilon^{\alpha_1+\alpha_2}}\right).
\ee
Therefore, we need to find combination $\{j_1, j_2\}$ which is the solution to the following optimization problem
\be
2j_1s_1 + 2j_2s_2 \xrightarrow{\min} j_1(2s_1+2\nu +\alpha_1)+j_2(2s_2 + \alpha_2) \geq \tau_{\varepsilon}, \ j_1, j_2 \geq 0.
\ee
It is easy to check that the solution is $\{j_1, j_2\}=\{\frac{\tau_{\varepsilon}}{2s_1 + 2\nu + \alpha_1}, 0 \}$, if $\alpha_2s_1 < s_2(2\nu + \alpha_1)$, and $\{j_1, j_2\}=\{0, \frac{\tau_{\varepsilon}}{2s_2 + \alpha_2}\}$, if $\alpha_2s_1 \geq s_2(2\nu + \alpha_1)$. Hence, the lower bounds are
\be
\label{delta1}
\delta^p =C A^p \left\{ \begin{array}{ll}
\left[\frac{\varepsilon^{2\overline{\alpha}}}{A^{2}}\right]^{\frac{ps_1}{2s_1+2\nu +\alpha_1}}, & \mbox{if}\ \ \alpha_2s_1 \leq s_2 (2\nu +\alpha_1),\\
\left[\frac{\varepsilon^{2\overline{\alpha}}}{A^{2}}\right]^{\frac{ps_2}{2s_2 +\alpha_2}},&  \mbox{if}\ \  \alpha_2s_1 > s_2 (2\nu +\alpha_1).
\end{array} \right.
\ee
\underline{\bf  The sparse-dense case.}
Using the same test functions $\tilde{f}_{j_1j_2}$ and $f_{j_1j_2}$, as in Benhaddou et al.~(2013), with a finitely supported basis $\eta_{j_2k_2}$ replaced by a band-limited basis $\psi_{j_2k_2}$, and following the same procedure as in the dense-dense case,  it can be shown that the lower bounds are 
\be \label{delta2}
\delta^p = C A^p\left\{ \begin{array}{ll}
\left[\frac{\varepsilon^{2\overline{\alpha}}}{A^{2}}\right]^{\frac{ps_2}{2s_2 +\alpha_2}},&  \mbox{if}\ \  \alpha_2s^*_1 \geq s_2 (2\nu+ \alpha_1-1) ,\\
\left[\frac{\varepsilon^{2\overline{\alpha}}}{A^{2}}\right]^{\frac{p(s_1+1/p-1/\pi)}{2s^*_1+2\nu +\alpha_1 -1}},&  \mbox{if}\ \  \ \alpha_2s^*_1 < s_2 (2\nu+ \alpha_1-1).
\end{array} \right.
\ee
To complete the proof, notice that the highest of the lower bounds corresponds to 
\be
d=\min \left\{ \frac{ps_1}{2s_1 + 2\nu +\alpha_1}, \frac{ps_2}{2s_2 +\alpha_2}, \frac{p(s_1+\frac{1}{p}-\frac{1}{\pi})}{2s^*_1+2\nu + \alpha_1-1} \right\}.
\ee
$\Box$\\
{\bf Proof of Theorem \ref{th:lowerbds-r}}.\\ \underline{\bf  The dense-dense case.}
Let $\eta$ be the matrix with components $\eta_{{\bf k}}=\{0, 1\}$, ${\bf k}=(k_1, k_2, \cdots, k_r)$, $k_i=0, 1, \cdots, 2^{j_i}-1$, $i=1, 2, \cdots, r$, and denote the set of all possible values $\eta$ by $H$. Define the functions
\be 
f_{{\bf j}}(t, x)=\gamma_{{\bf j}}\sum^{2^{j_1}-1}_{k_1=0}\sum^{2^{j_2}-1}_{k_2=0}\cdots \sum^{2^{j_r}-1}_{k_r=1}\eta_{{\bf k}}\psi_{j_1, k_1}(t)\Pi^{r-1}_{i=1}\psi_{j_i, k_i}(x_i).\label{tesf1-r}
\ee
Note that here $\card(H)=2^N$, with $N=2^{j_1+ j_2+\cdots + j_r}$. Hence, following the same steps as in the proof of Theorem \ref{th:lowerbds}, it can be shown that Note that the lower bounds are 
\be
\label{delta1-r}
\delta^p =C A^p \left\{ \begin{array}{ll}
\left[\frac{\varepsilon^{2\overline{\alpha}}}{A^{2}}\right]^{\frac{ps_1}{2s_1+2\nu +\alpha_1}}, & \mbox{if}\ \ \alpha_{2o}s_1 \leq s_{2o} (2\nu +\alpha_1),\\
\left[\frac{\varepsilon^{2\overline{\alpha}}}{A^{2}}\right]^{\frac{ps_{i}}{2s_{i} +\alpha_{i}}},&  \mbox{if}\ \  \alpha_{i}s_1 > s_{i} (2\nu +\alpha_1), \ \ \frac{s_i}{\alpha_i} < \min_{i\neq k}\{\frac{s_k}{\alpha_k}, k=2, 3, \cdots, r\},
\end{array} \right.
\ee
and $i=2, 3, \cdots, r$.\\
\underline{\bf  The sparse-dense case.} We apply similar approach to test functions 
\be 
f_{{\bf j}}(t, x)=\gamma_{{\bf j}}\sum^{2^{j_2}-1}_{k_2=0}\cdots \sum^{2^{j_r}-1}_{k_r=1}\eta_{{\bf k}}\psi_{j_1, k_1}(t)\Pi^{r-1}_{i=1}\psi_{j_i, k_i}(x_i).\label{tesf2-r}
\ee
Here $N=2^{ j_2+\cdots + j_r}$ and $\card(H)=2^N$ and it can be shown that the lower bounds are
\be \label{delta2-r}
\delta^p = C A^p\left\{ \begin{array}{ll}
\left[\frac{\varepsilon^{2\overline{\alpha}}}{A^{2}}\right]^{\frac{ps_i}{2s_i +\alpha_i}},&  \mbox{if}\ \  \alpha_is^*_1 \geq s_i (2\nu+ \alpha_1-1) ,\ \ \frac{s_i}{\alpha_i} < \min_{i\neq k}\{\frac{s_k}{\alpha_k}, k=2, 3, \cdots, r\}\\
\left[\frac{\varepsilon^{2\overline{\alpha}}}{A^{2}}\right]^{\frac{p(s_1+1/p-1/\pi)}{2s^*_1+2\nu +\alpha_1 -1}},&  \mbox{if}\ \  \ \alpha_is^*_1 < s_i (2\nu+ \alpha_1-1),
\end{array} \right.
\ee
and $i=2, 3, \cdots, r$. Hence, the lower bounds of the $L^p$-risk correspond to 
\be
d'=\min \left\{ \frac{ps_1}{2s_1 + 2\nu +\alpha_1},  \frac{p(s_1+\frac{1}{p}-\frac{1}{\pi})}{2s^*_1+2\nu + \alpha_1-1},\frac{ps_{2o}}{2s_{2o} +\alpha_{2o}}\right\},
\ee
and the pair $(s_{2o}, \alpha_{2o})$ is defined in \fr{s-op}. $\Box$
\subsection{Proof of the upper bounds.}
{\bf Proof of Lemma \ref{lem:Var}.} Define the quantities
\be   \label{devbet}
\tilde{\beta}_{\omega}-\beta_{\omega}=\varepsilon^{\overline{\alpha}} \sum_{m_1\in W_{j_1}}\sum_{m_2\in W_{j_2}}\overline{\tilde{\psi}_{j_1k_1}(m_1)}\tilde{\psi}_{j_2k_2}(m_2) \frac{\tilde{Z}^{\alpha}(m_1, m_2)}{\tilde{g}(m_1, m_2)}.
\ee 
Consider the Riesz poly-potential operator
\be \label{rpoten}
I^{\alpha}f(t, x)=\frac{1}{\gamma_2(\alpha)}\int_{R^2}\frac{f(t_2, x_2)}{|t_1-t_2|^{1/2+\alpha_1/2}|x_1-x_2|^{1/2+\alpha_2/2}}dt_2dx_2,
\ee
where $\gamma_2(\alpha)=4\Gamma\left(\frac{1}{2}-\frac{\alpha_1}{2}\right) \Gamma\left(\frac{1}{2}-\frac{\alpha_2}{2}\right)cos((1-\alpha_1)\frac{\pi}{4}) cos((1-\alpha_2)\frac{\pi}{4})$. Then, the anisotropic two-dimensional fractional Brownian sheet $Z^{\alpha}(t, x)=dW^{\alpha}(t, x)$ allows the wavelet-vaguelette representation 
\be \label{wavag}
dW^{\alpha}(t, x)=\sum_{j_1, k_1}\sum_{j_2, k_2}\zeta_{j_1k_1;j_2k_2} I^{\alpha}\left\{\psi_{j_1k_1}(t)\psi_{j_2k_2}(x)\right\},
\ee 
where $\zeta_{j_1k_1;j_2k_2}$ are white noise processes, and $\psi_{jk}$ is a Meyer-type wavelet basis. Then applying the two dimensional Fourier transform to \fr{wavag} yields
\be \label{fourwavag}
\tilde{W^{\alpha}}({\bf{m}})= |m_1|^{\frac{\alpha_1-1}{2}}|m_2|^{\frac{\alpha_2-1}{2}}\sum_{j_1,k_1}\sum_{j_2,k_2}\zeta_{j_1k_1;j_2k_2}\overline{\tilde{\psi}_{j_1k_1}(m_1)}\tilde{\psi}_{j_2k_2}(m_2),
\ee
where ${\bf{m}}=(m_1, m_2)$. Let us evaluate the covariance of \fr{fourwavag}. Indeed,
\be
\Cov[\tilde{Z^{\alpha}}({\bf{m}}), \tilde{Z^{\alpha}}({\bf{l}})]= |m_1l_1|^{\frac{\alpha_1-1}{2}}|m_2l_2|^{\frac{\alpha_2-1}{2}}\sum_{\omega\in\Omega}\overline{\tilde{\psi}_{j_1k_1}(m_1)}\tilde{\psi}_{j_2k_2}(m_2)\overline{\tilde{\psi}_{j_1k_1}(l_1)}\tilde{\psi}_{j_2k_2}(l_2).
\ee 
Consequently,
\be \label{covf}
\Cov[\tilde{Z^{\alpha}}({\bf{m}}), \tilde{Z^{\alpha}}({\bf{l}})]= |m_1l_1|^{\frac{\alpha_1-1}{2}}|m_2l_2|^{\frac{\alpha_2-1}{2}}\left(\sum_{j_1;k_1} \tilde{\psi}_{j_1k_1}(m_1) \tilde{\psi}_{j_1k_1}(l_1)\right) \left(\sum_{j_2;k_2} \tilde{\psi}_{j_2k_2}(m_2) \tilde{\psi}_{j_2k_2}(l_2)\right).
\ee 
Now, evaluating the magnitude of \fr{covf}, and using H{\"o}lder's Inequality along with the fact that $|\tilde{\psi}_{jk}(m)| \leq 2^{-j/2}$, yields
\beqns
| \Cov[\tilde{Z^{\alpha}}({\bf{m}}), \tilde{Z^{\alpha}}({\bf{l}})]|^2&\leq& |m_1l_1|^{{\alpha_1-1}}|m_2l_2|^{{\alpha_2-1}}\left(\sum_{j_1;k_1} |\tilde{\psi}_{j_1k_1}(m_1)|^2 |\tilde{\psi}_{j_1k_1}(l_1)|^2\right) \left(\sum_{j_2;k_2} |\tilde{\psi}_{j_2k_2}(m_2)|^2 |\tilde{\psi}_{j_2k_2}(l_2)|^2\right)\nonumber\\
& \leq& 4|m_1l_1|^{{\alpha_1-1}}|m_2l_2|^{{\alpha_2-1}}.
\eeqns 
Hence,
\be \label{covab}
| \Cov[\tilde{Z^{\alpha}}({\bf{m}}), \tilde{Z^{\alpha}}({\bf{l}})]|\leq 2 |m_1l_1|^{\frac{\alpha_1-1}{2}}|m_2l_2|^{\frac{\alpha_2-1}{2}}.
\ee
Finally, let us evaluate the variance of $\tilde{\beta}_{\omega}$. Indeed, using \fr{devbet} and \fr{covab}, the variance has the form 
\beqn
\EE\left|\tilde{\beta}_{\omega}-\beta_{\omega}\right|^2&=& \varepsilon^{\alpha_1+\alpha_2}\sum_{m_1, m_2}\sum_{l_1, l_2}\overline{\tilde{\psi}_{j_1k_1}(m_1)}\tilde{\psi}_{j_2k_2}(m_2) \overline{\tilde{\psi}_{j_1k_1}(l_1)}\tilde{\psi}_{j_2k_2}(l_2)\frac{\EE[\overline{\tilde{Z^{\alpha}}(m_1, m_2)}\tilde{Z^{\alpha}}(l_1, l_2)]}{\overline{\tilde{g}(m_1, m_2)}\tilde{g}(l_1, l_2)}\nonumber\\
&\leq& 2\varepsilon^{\alpha_1+\alpha_2}\left(\sum_{m_1, m_2}{\tilde{\Psi}_{\omega}(m_1, m_2)}\frac{|m_1|^{\alpha_1-1}|m_2|^{\alpha_2-1}}{\tilde{g}(m_1, m_2)}\right) \left(\sum_{l_1, l_2}{\tilde{\Psi}_{\omega}(l_1, l_2)}\frac{|l_1|^{\alpha_1-1}|l_2|^{\alpha_2-1}}{\tilde{g}(l_1, l_2)}\right)\nonumber\\
&=& 2\varepsilon^{\alpha_1+\alpha_2}\left(\sum_{m_1, m_2}\overline{\tilde{\psi}_{j_1k_1}(m_1)}\tilde{\psi}_{j_2k_2}(m_2)\frac{|m_1|^{\alpha_1-1}|m_2|^{\alpha_2-1}}{\tilde{g}(m_1, m_2)}\right)^2\nonumber\\
&\leq& 2\varepsilon^{\alpha_1+\alpha_2}C^{-2}_1\left(\frac{8\pi}{3}\right)^{2\nu +\alpha_1+\alpha_2-2}2^{j_1(2\nu +\alpha_1-1)+j_2(\alpha_2-1)}. \label{varproo} 
\eeqn
Results \fr{var} and  \fr{varp}  follow from properties of Gaussian random variables. $\Box$\\
{\bf Proof of Lemma \ref{lemmabetap}.}
First, note that under {\bf Assumption A.3.}, one has
\begin{equation}\label{eq: 3.19}
\sum\limits_{k_1,k_2} \mid \beta_{j_1,k_1,j_2,k_2} \mid^{\pi} \leq A^{\pi} 2^{-\pi[(j_1s_1+j_2s_2)+ (\frac{1}{2} - \frac{1}{\pi} )(j_1+j_2)]}, \ \forall j_1, j_2\geq0.
\end{equation}
If $p\geq \pi$, one has 
\beqns
\sum\limits_{k_1,k_2} \mid \beta_{j_1,k_1,j_2,k_2} \mid^p 
&\leq&\sum\limits_{k_1,k_2} \left(    \mid \beta_{j_1,k_1,j_2,k_2} \mid^{\pi} \left\{ \max |\beta_{j_1,k_1,j_2,k_2} |^{\pi} \right\}^{\frac{p-\pi}{\pi}} \right) \nonumber\\
&=&\left\{ \max |\beta_{j_1,k_1,j_2,k_2} |^{\pi} \right\}^{\frac{p-\pi}{\pi}} \sum\limits_{k_1,k_2}   \mid \beta_{j_1,k_1,j_2,k_2} \mid^{\pi}  \nonumber\\
&\leq&A^p 2^{- p[(j_1s_1+j_2s_2)+ (\frac{1}{2} - \frac{1}{\pi} )(j_1+j_2) ]}.
\eeqns
If $p<\pi$, then by H{\"o}lder's Inequality, one has
\beqns
\sum\limits_{k_1,k_2} \mid \beta_{j_1,k_1,j_2,k_2} \mid^p &\leq& \left(   \sum\limits_{k_1,k_2} \mid \beta_{j_1,k_1,j_2,k_2} \mid^{\pi}  \right) ^{\frac{p}{\pi}}   \left( \sum\limits_{k_1,k_2} 1\right)^{1- \frac{p}{\pi}} \nonumber\\
&\leq& A^p 2^{-p[(j_1s_1+j_2s_2)+ (\frac{1}{2} - \frac{1}{\pi} )(j_1+j_2) ]}  2^{(j_1+j_2)(1-\frac{p}{\pi})} \nonumber\\
&\leq&A^p 2^{- p[(j_1s_1+j_2s_2)+ (\frac{1}{2} - \frac{1}{p} )(j_1+j_2) ]}.
\eeqns
Combining the above results completes the proof. $\Box$\\
{\bf Proof of Lemma \ref{lem:Lar-D}.} Since the quantities \fr{devbet} are zero mean Gaussian random variables having variance of form \fr{varproo}, by the Gaussian tail probability inequality, one has 
\beqn
\Pr(\Theta_{\omega, \eta})&\leq& 2\Pr \left(Z> \eta \gamma \sqrt{|\ln(\varepsilon)|}\left(\frac{C^2_1}{2}\right)^{1/2}\left(\frac{3}{8\pi}\right)^{\nu +\alpha_1/2 +\alpha_2/2-1}\right)\nonumber\\
&=& \frac{C}{\sqrt{|\ln(\varepsilon)|}}\exp\left\{-2\overline{\alpha}|\ln(\varepsilon)|\tau\right\},
\eeqn
where $Z$ is a standard normal, $C$ is some explicit positive constant and $\tau$ appears in Lemma 2. $\Box$\\
{\bf Proof of Theorem \ref{th:upperbds}.} The proof follows a procedure which is similar to that in  Benhaddou et al.~(2013) with the adaptation to the  $L^p$-norm case. Denote
\be  \label{Lev:chi}
\chi_{\varepsilon}= \left[{A^{-2}} {\varepsilon^{2\overline{\alpha}}}|\ln(\varepsilon)|\right],\ \
2^{j_{i0}}= \left[\chi_{\varepsilon}\right]^{-\frac{d}{ps"_i}},\ \
i=1, 2,
\ee
and observe that with the choice of $J_1$ and $J_2$ given by \fr{Lev:J}, the $L^p$-risk allows the decomposition
\be \label{eseroll}
\EE \| \widehat{f}_{\epsilon}-f\|^p_p\leq R_1+ R_2 + R_3 + R_4,
\ee
where 
\beqns
R_1&=& 2^{2(p-1)}\EE \left\|\sum^{2^{m_{0}}-1}_{k_1=0}\sum^{2^{m_{0}}-1}_{k_2=0} (\tilde{ \beta}_{m_{0}, k_1, m_{0}, k_2}-\beta_{m_{0}, k_1, m_{0}, k_2})\psi_{m_{0}, k_1}(t)\psi_{m_{0}, k_2}(x)\right\|_p^{p},\nonumber\\
R_2&=& 2^{3(p-1)} \EE \int \int \left(\sum_{\bom \in \Om(J_1, J_2)} \left| \tilde{ \beta}_{\bom}-\beta_{\bom}\right|^2 \II \left(  \left| \tilde{ \beta}_{\bom}  \right| >  \lambda^{\alpha}_{j;\varepsilon} \right)\psi^2_{j_1, k_1}(t)\psi^2_{j_2, k_2}(x)\right)^{p/2}dtdx, \nonumber\\
R_3&=& 2^{3(p-1)} \EE \int \int \left(\sum_{\bom \in \Om(J_1, J_2)} \left| { \beta}_{\bom}  \right|^2 \II \left( \left| \tilde{ \beta}_{\bom}  \right| < \lambda^{\alpha}_{j;\varepsilon} \right)\psi^2_{j_1, k_1}(t)\psi^2_{j_2, k_2}(x)\right)^{p/2}dtdx,\\
R_4&=& 2^{(p-1)}\int \int \left(  \left(\sum^{\infty}_{j_1=J_1}\sum^{\infty}_{j_2=m_{0}+1}+ \sum^{\infty}_{j_1=m_{0}+1}\sum^{\infty}_{j_2=J_2}\right)\sum^{2^{j_1}-1}_{k_1=0}\sum^{2^{j_2}-1}_{k_2=0} \left| { \beta}_{\bom}  \right|^2 \psi^2_{j_1, k_1}(t)\psi^2_{j_2, k_2}(x)  \right)^{p/2}dtdx.
\eeqns
For $R_1$ and $R_4$, using \fr{var} in $R_1$ and \fr{eq11} in $R_4$, as $\varepsilon \rightarrow 0$, yields
\be  \label{R_1}
R_1+R_4=O \left( \varepsilon^{2\overline{\alpha}}+ \sum^{\infty}_{j_1=J_1}A^p2^{-pj_1s''_1}+ \sum^{\infty}_{j_2=J_2}A^p2^{-pj_2s''_2}\right)=O \left( A^p \left[ \chi_{\varepsilon}\right]^d \right).
\ee
Using Minkowski's Inequality for $p \geq 2$ and Jensen's Inequality for $p<2$, and Temlyakov property, $R_2$ and $R_3$ can be partitioned as $R_2\leq R_{21} + R_{22}$ and  $R_3 \leq R_{31} + R_{32}$, where
\beqn  
R_{21}&=& O \left( \sum_{\bom \in \Om(J_1, J_2)}2^{(j_1+j_2)(\frac{p}{2}-1)}\EE\left| \tilde{ \beta}_{\bom}-\beta_{\bom}\right|^p \left[\Pr \left(  \left| \tilde{ \beta}_{\bom}-\beta_{\bom}  \right| > \frac{1}{2} \lambda^{\alpha}_{j;\varepsilon}\right)\right]^{1/2}  \right),\label{r21}\\
R_{22}&=& O \left(\sum_{\bom \in \Om(J_1, J_2)}2^{(j_1+j_2)(\frac{p}{2}-1)}\EE \left| \tilde{ \beta}_{\bom}-\beta_{\bom}\right|^p \II \left(  \left|  \beta_{\bom} \right| >  \frac{1}{2}\lambda^{\alpha}_{j;\varepsilon} \right) \right),\label{r22}\\
R_{31}&=& O \left(\sum_{\bom \in \Om(J_1, J_2)}2^{(j_1+j_2)(\frac{p}{2}-1)}\left| { \beta}_{\bom}  \right|^p \Pr \left( \left| \tilde{ \beta}_{\bom}-\beta_{\bom}  \right| > \frac{1}{2} \lambda^{\alpha}_{j;\varepsilon} \right)\right),\label{r31}\\
R_{32}&=& O \left(\sum_{\bom \in \Om(J_1, J_2)}2^{(j_1+j_2)(\frac{p}{2}-1)} \left| { \beta}_{\bom}  \right|^p\II \left(  \left|  \beta_{\bom} \right| <  \frac{3}{2}\lambda^{\alpha}_{j;\varepsilon} \right)\right).\label{r32}
\eeqn
Combining \fr{r21} and \fr{r31}, similar calculations as in Benhaddou et al.~(2013), yield
\be
R_{21} + R_{31} = O  \left( A^p\left( {\varepsilon^{2\overline{\alpha}}}\right)^{\frac{\tau}{2}-\frac{p}{2}}  \right)=  O \left( A^p \left[ \chi_{\varepsilon}\right]^d \right). \label{r21r31}
\ee
Now, combining \fr{r22} and \fr{r32}, and using \fr{var} and \fr{Thres}, gives
\be  \label{r22r32}
 \Delta=R_{22} + R_{32}=  O \left( \sum_{\bom \in \Om(J_1, J_2)}2^{(j_1+j_2)(\frac{p}{2}-1)}\min \left\{ \left| { \beta}_{\bom}  \right|^p,   2^{\frac{p}{2}[j_1(2\nu+\alpha_1-1)+j_2(\alpha_2-1)]} \left[ A^{-2} {\varepsilon^{2\overline{\alpha}}}|\ln(\varepsilon)|\right]^{\frac{p}{2}}\right\} \right).
\ee
Finally, $\Delta$ can be decomposed into the following components
\beqn  
 \Delta_1&=&  O \left( \left\{\sum^{J_1-1}_{j_1=j_{10}+1}\sum^{J_2-1}_{j_2=m_{0}}+ \sum^{J_1-1}_{j_1=m_{0}}\sum^{J_2-1}_{j_2=j_{20}+1}\right\}A^p 2^{-j_1ps''_1} 2^{-j_2ps''_2} \right), \label{del1}\\
 \Delta_2&=&O \left( \sum^{j_{10}}_{j_1=m_{0}} \sum^{j_{20}}_{j_2=m_{0}}A^p  2^{\frac{p}{2}(j_1(2\nu +\alpha_1)+\alpha_2j_2)}\left[ A^{-2} {\varepsilon^{2}}|\ln(\varepsilon)|\right]^{\frac{p}{2}} \II \left(\Xi \right) \right), \label{del2}\\
 \Delta_3&=& O \left( \sum^{j_{10}}_{j_1=m_{0}} \sum^{j_{20}}_{j_2=m_{0}}A^{p'}2^{-p'(j_1s''_1+j_2s''_2)} \left[ A^p  [\chi_{\varepsilon}]^{\frac{p}{2}} 2^{j_1(\frac{p}{2}(2\nu +\alpha_1) -1)+j_2(\alpha_2\frac{p}{2}-1)} \right] ^{1-\frac{p'}{p}}\II \left(\Xi^c \right)\right)\label{del3}
\eeqn
where $\Xi=\left\{ j_1, j_2: 2^{\frac{p}{2}(j_1(2\nu +\alpha_1)+\alpha_2j_2)} \leq   \left[ \chi_{\varepsilon}\right]^{d-\frac{p}{2}} \right \}$.\\
 \underline{Case 1: $  s_1  \geq \frac{s_2}{\alpha_2}(2\nu+\alpha_1)$  and  $\frac{s_2}{\alpha_2} > \frac{p}{2}(\frac{1}{p'}-\frac{1}{p}) $.} In this case, $d=\frac{p s_2}{2s_2+\alpha_2}  $, and
\beqn
			\Delta_3
			&=&O \left(   A^p [\chi_{\varepsilon}]^{\frac{p}{2}(1-\frac{p'}{p})}\sum\limits_{j_1=m_{0}}^{j_{10}} 2^{-p' j_1[s_1-\frac{p}{2}(\frac{1}{p'}-\frac{1}{p})(2\nu + \alpha_1)]}  \sum\limits_{j_2=m_{0}}^{j_{20}} 2^{-j_2p'[s_2-\frac{p}{2}\alpha_2(\frac{1}{p'}-\frac{1}{p}) ]} \II\left(  \Xi^c \right)  \right) \nonumber\\
			&=&O \left(  A^p    [\chi_{\varepsilon}]^{\frac{p s_2}{2s_2+\alpha_2} } \sum\limits_{j_1=m_{0}}^{j_{10}} 2^{  {-p'j_1}[s_1- \frac{s_2}{\alpha_2}(2\nu+\alpha_1)] } \right)\nonumber\\
			&=&O\left( A^p [\chi_{\varepsilon}]^{\frac{p s_2}{2s_2+\alpha_2} }    \left[\left| \ln(\epsilon)  \right|\right]^{\II(  s_1 = \frac{s_2}{\alpha_2}(2\nu+\alpha_1) )} \right).
\eeqn
\underline{Case 2:  $(2\nu+\alpha_1)\frac{p}{2}(\frac{1}{\pi}-\frac{1}{p})< s_1< \frac{s_2}{\alpha_2}(2\nu+\alpha_1)$}. In this case, $d=\frac{p s_1}{2s_1+2\nu+\alpha_1} $, and
\beqn
	\Delta_3&=&O \left(   A^p [\chi_{\varepsilon}]^{\frac{p}{2} (1-\frac{p'}{p})}\sum\limits_{j_1=m_{0}}^{j_{10}} 2^{ -p' j_1[ s_1-(2\nu+\alpha_1)\frac{p}{2}(\frac{1}{p'}-\frac{1}{p}) ]} \sum\limits_{j_2=m_{0}}^{j_{20}}2^{-j_2p'[s_2-\frac{p}{2}\alpha_2(\frac{1}{p'}-\frac{1}{p}) ]} \II\left(  \Xi^c \right)  \right) \nonumber\\
		         &=&O \left(   A^p    [\chi_{\varepsilon}]^{\frac{p s_1}{2s_1+2\nu+\alpha_1}} \sum\limits_{j_2=m_{0}}^{j_{20}}2^{-j_2p' \left(s_2- \frac{\alpha_2s_1}{2\nu+\alpha_1}  \right)   } \right)\nonumber\\
		         &=&O\left(A^p    [\chi_{\varepsilon}]^d \right).
\eeqn
\underline{Case 3:  $s_1\leq   (2\nu+\alpha_1)\frac{p}{2}(\frac{1}{\pi}-\frac{1}{p})$ and $\frac{s_2}{\alpha_2} > \frac{p}{2}(\frac{1}{p'}-\frac{1}{p})$}. In this case, $d= \frac{p(s_1 + \frac{1}{p} -\frac{1}{\pi})}{2s^*_1+2\nu +\alpha_1-1}$,  and
\beqn 
			\Delta_3&=&O \left(  A^p [\chi_{\varepsilon}]^{\frac{p}{2} (1-\frac{p'}{p})}\sum\limits_{j_1=m_{0}}^{j_{10}} 2^{-p' j_1[s_1-\frac{p}{2}(\frac{1}{p'}-\frac{1}{p})(2\nu + \alpha_1)]} \sum\limits_{j_2=m_{0}}^{j_{20}} 2^{-p' j_2[s_2-\frac{p}{2}(\frac{1}{p'}-\frac{1}{p})\alpha_2]} 			\right)\nonumber\\
			&=&O \left(  A^p [\chi_{\varepsilon}]^{\frac{p(s_1 + \frac{1}{p} -\frac{1}{\pi})}{2s^*_1+2\nu +\alpha_1-1}} \left[\left| \ln(\epsilon)\right|\right]^{\II \left(s_1= (2\nu+\alpha_1)\frac{p}{2}(\frac{1}{p'}-\frac{1}{p}) \right)+\II \left(s_2= \alpha_2\frac{p}{2}(\frac{1}{p'}-\frac{1}{p}) \right)} \right). \label{del3b}
			\eeqn
		Combining the results from \fr{R_1} to \fr{del3b} completes the proof. $\Box$ \\
{\bf Proof of Lemma \ref{lem:Var_r}.}
The proof will be similar to the proof of {\bf Lemma 1}. Define the quantities
\be   \label{devbet-r}
\tilde{\beta}_{\omega_r}-\beta_{\omega_r}=\varepsilon^{\overline{\alpha}} \sum_{m_1\in W_{j_1}}\cdots \sum_{m_r\in W_{j_r}}\tilde{\Psi}_{\omega_r}({\bf m_r}) \frac{\tilde{Z}^{\alpha}({\bf m_r})}{\tilde{g}({\bf m_r})}.
\ee 
Consider the Riesz poly-potential operator
\be \label{rpoten-r}
I^{\alpha}f(t, x_1,\ldots, x_{r-1})=\frac{1}{\gamma_r(\alpha)}\int_{R^{r}}\frac{f(t', x'_1,\ldots, x'_{r-1}  )}{|t-t'|^{1/2+\alpha_1/2}  \Pi^{r}_{i=2}|x_{i-1}-x'_{i-1}|^{1/2+\alpha_i/2}   }dt'dx'_1\cdots dx'_{r-1},
\ee
where $\gamma_r(\alpha)=2^{r}\Pi^r_{i=1}\left[\Gamma\left(\frac{1}{2}-\frac{\alpha_i}{2}\right)cos((1-\alpha_i)\frac{\pi}{4}) \right]$. Then, the anisotropic $r$-dimensional fractional Brownian sheet $Z^{\alpha}(t, x_1,\ldots, x_{r-1})=dW^{\alpha}(t, x_1,\ldots, x_{r-1})$ allows the wavelet-vaguelette representation 
\be \label{wavag-r}
dW^{\alpha}(t, x_1,\ldots, x_{r-1})=\sum_{\omega_r \in \Omega_r}\zeta_{\omega_r} I^{\alpha}\left\{\Psi_{\omega_r}(t, {\bf x})\right\},
\ee 
Applying the two dimensional Fourier transform to \fr{wavag-r}, yields
\be \label{fourwavag-r}
\tilde{W^{\alpha}}({\bf m_r})= \Pi^r_{i=1}|m_i|^{\frac{\alpha_i-1}{2}} \sum_{\omega_r \in \Omega_r}\zeta_{\omega_r} \tilde{\Psi}_{\omega_r}({\bf m_r}),
\ee
The covariance of \fr{fourwavag-r} is given by
\be \label{covf-r}
\Cov[\tilde{Z^{\alpha}}({\bf m_r}), \tilde{Z^{\alpha}}({\bf l_r })]=\Pi^r_{i=1} \left(|m_il_i|^{\frac{\alpha_i-1}{2}}\sum_{j_i;k_i} \tilde{\psi}_{j_ik_i}(m_i) \tilde{\psi}_{j_ik_i}(l_i)\right).
\ee 
Evaluating the magnitude of \fr{covf-r}, yields
\be \label{covab-r}
| \Cov[\tilde{Z^{\alpha}}({\bf m_r}), \tilde{Z^{\alpha}}({\bf l_r})]|\leq 2^{\frac{r}{2}} \Pi^r_{i=1} |m_il_i|^{ \frac{\alpha_i-1}{2}}.
\ee
Using \fr{devbet-r} and \fr{covab-r}, similar to \fr{varproo}, the variance has the form 
\beqn
\EE\left|\tilde{\beta}_{\omega_r}-\beta_{\omega_r}\right|^2&=& \varepsilon^{2 \overline{\alpha}} \sum_{m_1,\ldots, m_r}\sum_{l_1,\ldots, l_r}\overline{\tilde{\Psi}_{\omega_r}({\bf m_r})} \tilde{\Psi}_{\omega_r}({\bf l_r})\frac{\EE[\overline{\tilde{Z^{\alpha}}({\bf m_r})}\tilde{Z^{\alpha}}({\bf l_r})]}{\overline{\tilde{g}({\bf m_r})}\tilde{g}({\bf l_r})}\nonumber\\
&\leq& K \varepsilon^{2 \overline{\alpha}} 2^{j_1(2\nu +\alpha_1-1)+j_2(\alpha_2-1) +\cdots+ j_r(\alpha_r-1) }. \label{varproo-r} 
\eeqn
Finally, results \fr{var_r} and  \fr{varp_r}  follow from properties of Gaussian random variables. $\Box$\\
{\bf Proof of Lemma \ref{lem:Lar-D-r}.} The proof is very similar to that of Lemma \ref{lem:Lar-D}, so we skip it. $\Box$

\end{document}